\documentclass{amsart}
\usepackage[final]{epsfig}
\usepackage{graphics}
\usepackage{float}
\usepackage{amsmath}
\usepackage{amsfonts}
\usepackage{latexsym}
\usepackage{amssymb}
\usepackage{graphicx}
\usepackage{url}
\usepackage{epstopdf}
\usepackage{verbatim}
\usepackage[margin=1in]{geometry}
\usepackage{amsthm}
\usepackage{enumitem}
\usepackage{tikz}
\usepackage[stable]{footmisc}
\usepackage{bm}
\usepackage{tabulary}
\usepackage{booktabs}

\makeatletter
\newtheorem*{rep@theorem}{\rep@title}
\newcommand{\newreptheorem}[2]{%
\newenvironment{rep#1}[1]{%
 \def\rep@title{#2 \ref{##1}}%
 \begin{rep@theorem}}%
 {\end{rep@theorem}}}
\makeatother

\newtheorem{lemma}{Lemma}[section]
\newtheorem{proposition}[lemma]{Proposition}

\newtheorem{remark}[lemma]{Remark}
\newtheorem{example}[lemma]{Example}
\newtheorem{theorem}[lemma]{Theorem}
\newtheorem{definition}[lemma]{Definition}
\newtheorem*{definition*}{Definition}
\newtheorem{corollary}[lemma]{Corollary}

\newtheorem*{theorem*}{Theorem}
\newreptheorem{theorem}{Theorem}
\newreptheorem{lemma}{Lemma}
\newreptheorem{definition}{Definition}

\newcommand{\proofend}{$\Box$\bigskip}

\newcommand{\R}{{\mathbb R}}

\def\proof{\paragraph{Proof.}}

\newcommand{\thmrestatement}[1] {\vspace{2 mm} \noindent \textbf{Restatement of Theorem #1. }\vspace{2 mm}}

\newcommand{\defrestatement}[1] {\vspace{2 mm} \noindent \textbf{Restatement of Definition #1. }\vspace{2 mm}}

\begin{document}

\title{Tropical Spectral Theory of Tensors}

\author{Emmanuel Tsukerman
}

\date{\today}
\address{Department of Mathematics, University of California,
Berkeley, CA 94720-3840}
\email{e.tsukerman@berkeley.edu}

\begin{abstract}
We introduce and study tropical eigenpairs of tensors, a generalization of the tropical spectral theory of matrices. We show the existence and uniqueness of an eigenvalue. We associate to a tensor a directed hypergraph and define a new type of cycle on a hypergraph, which we call an H-cycle. The eigenvalue of a tensor turns out to be equal to the minimal normalized weighted length of H-cycles of the associated hypergraph. We show that the eigenvalue can be computed efficiently via a linear program. Finally, we suggest possible directions of research.
\end{abstract}

\maketitle
\setcounter{tocdepth}{1}

\section{Background}

A tensor of order $m$ and rank $n$ is an array $A=(a_{i_1 \cdots i_m})$ of elements of a field $K$ (which we shall take to be $\R$ or $\overline{\R}=\R \cup \{\infty\}$), where $1 \leq i_1,\ldots,i_m \leq n$. In ordinary arithmetic, given $x \in \R^n$, we define
\[
(Ax^{m-1})_i:=\sum_{i_2,\ldots,i_m=1}^n a_{i i_2 \cdots i_m} x_{i_2} \cdots x_{i_m}.
\] 

An H-eigenpair \cite{MR2178089} of a tensor is defined as follows. Define $x^{[m-1]}=(x_i^{m-1})_{i}$. Then an H-eigenpair is a pair $(x,\lambda) \in \mathbb{P}^{n-1} \times \R$ such that
\[
Ax^{m-1}=\lambda x^{[m-1]}.
\]

Let $A$ be a $n \times n$ matrix with entries in the tropical semiring $(\R,\oplus,\odot)$. An eigenvalue of $A$ is a number $\lambda$ such that
\[
A \odot \bm{v}=\lambda \odot \bm{v}.
\]
The nature of tropical eigenpairs is understood in the setting of matrices (\cite{MR3033951},\cite{MR3201242}) but a survey of the literature shows no prior research on tropical eigenpairs of tensors.

\begin{definition}\normalfont
A \emph{tropical H-eigenpair} for a tensor $(a_{i_1 \cdots i_m}) \in \R^{n^m}$ of order $m$ and rank $n$ is a pair $(x,\lambda) \in \R^{n} / \R(1,1,\ldots,1) \times \R$ such that
\begin{align} \label{eigenProblemH}
\bigoplus_{i_2,\ldots,i_m=1}^n a_{i i_2 \cdots i_m} \odot x_{i_2} \odot \cdots \odot x_{i_m}=\lambda \odot x_i^{m-1}, \quad i=1,2,\ldots,n.
\end{align}
We call $x$ a \emph{tropical H-eigenvector} and $\lambda$ a \emph{tropical H-eigenvalue}.
\end{definition}

In the classical setting, several other definitions of eigenpairs of tensors exist. For instance, an  E-eigenpair is defined via the condition 
\[
Ax^{m-1}=\lambda x.
\]

We define the tropicalization here in an analogous manner. In this paper, we focus on H-eigenpairs and only discuss E-eigenpairs for purposes of comparison.

We take $\oplus$ to be min throughout.

\begin{example}
Take $n=2$ and $m=3$. Then a tropical H-eigenpair $(x,\lambda)$ satisfies
\begin{equation}
\begin{array}{l}
\min \{a_{111}+2x_1,a_{112}+x_1+x_2,a_{121}+x_2+x_1,a_{122}+2x_2\}=\lambda+2x_1 \\
\min \{a_{211}+2x_1,a_{212}+x_1+x_2,a_{221}+x_2+x_1,a_{222}+2x_2\}=\lambda+2x_2.
\end{array} 
\label{eq:d32}
\end{equation}
\end{example}

Without loss of generality, we will assume from now on that all tensors are symmetric in their last $m-1$ coordinates. That is,
\[
a_{i_1 i_2 \cdots i_m}=a_{i_1 i_{\sigma(2)} \cdots i_{\sigma(m)})}
\]
for all permutations $\sigma$ of $\{2,\ldots,m\}$. This is because only the smallest element in each orbit plays a role in the eigenproblem.

\begin{example}
Take $n=3$ and $m=3$. Then a tropical H-eigenpair $(x,\lambda)$ satisfies
\begin{equation}
\begin{array}{l}
 \min \{a_{111}+2x_1,a_{112}+x_1+x_2,a_{113}+x_1+x_3,a_{122}+2x_2,a_{123}+x_2+x_3,a_{133}+2x_3\}=\lambda+2x_1 \\
 \min \{a_{211}+2x_1,a_{212}+x_1+x_2,a_{213}+x_1+x_3,a_{222}+2x_2,a_{223}+x_2+x_3,a_{233}+2x_3\}=\lambda+2x_2 \\
 \min \{a_{311}+2x_1,a_{312}+x_1+x_2,a_{313}+x_1+x_3,a_{322}+2x_2,a_{323}+x_2+x_3,a_{333}+2x_3\}=\lambda+2x_3.
\end{array} 
\label{eq:n3m3}
\end{equation}
\end{example}

\section{Main Results}

We show that

\begin{theorem}\label{uniquenessExistence}
A tensor $A \in \R^{n^m}$ has a unique tropical H-eigenvalue $\lambda(A) \in \R$.
\end{theorem}

The result is all the more striking when compared with the situation of E-eigenpairs.  Experimentally, in 5000 runs of randomly generated symmetric $3 \times 3 \times 3$-tensors, we obtained the following distribution on the number of tropical E-eigenpairs  (from 0 eigenpairs to 7): [0, 4007, 7, 950, 0, 6, 1, 29].

Uniqueness and existence for tropical  H-eigenpairs can be extended to the case when $A \in \overline{\R}^{n^m}$ under suitable technical assumptions:
\begin{theorem}\label{uniquenessExistenceWithInfinity}
Let $A \in \overline{\R}^{n^m}$ and assume that for each $i=1,2,\ldots,n$, the sets
\[
S_i:=\{\{i_2^{(i)},\ldots,i_m^{(i)}\} \, : \, a_{i i_2^{(i)} \cdots i_m^{(i)}} \neq \infty\}
\]
are nonempty and mutually equal. Then $A$ has a unique tropical H-eigenvalue $\lambda(A) \in \R$.
\end{theorem}

Surprisingly, it turns  out that the H-eigenvalue is a solution to a much simpler problem given by a linear program which can be interpreted as merely requiring that $\lambda(A)$ be a subeigenvalue. This turns out to be a consequence of a special property of the H-cycle polytope (to be defined later).

\begin{theorem}\label{LPthm}
Assume the hypotheses of Theorem \ref{uniquenessExistenceWithInfinity}. The solution to the linear program
\begin{equation} \label{LP}
\begin{array}{ll@{}ll}
\text{maximize}  & \lambda &\\
\text{subject to}& a_{i_1 i_2 \cdots i_m}+x_{i_2}+\ldots+x_{i_m} \geq \lambda + (m-1)x_{i_1}, &   & \forall (i_1,i_2,\ldots,i_m) \in [n]^m\\
\end{array}
\end{equation}
  is equal to the H-eigenvalue of $A$. Dually, the H-eigenvalue of $A$ is given by
  \begin{equation} \label{LPdual}
\begin{array}{ll@{}ll}
\text{minimize}  & \sum_{(i_1,i_2,\ldots,i_m) \in [n]^m} a_{i_1 i_2 \cdots i_m} y_{i_1 i_2 \cdots i_m} &\\
\text{subject to}& \sum_{(i_1,i_2,\ldots,i_m) \in [n]^m} y_{i_1 i_2 \cdots i_m}((m-1)e_{i_1}-e_{i_2}-\ldots-e_{i_m})=0 \\
& \sum_{(i_1,i_2,\ldots,i_m) \in [n]^m} y_{i_1 i_2 \cdots i_m}=1 \\
& y_{i_1 i_2 \cdots i_m} \geq 0.
\end{array}
\end{equation}
\end{theorem}

For matrices ($m=2$), the dual problem has the interpretation of giving the minimal normalized length of any directed cycle in the weighted directed graph associated to the matrix $A$. Motivated by this, we associate to a tensor $A$ a weighted  $F$-hypergraph.

We summarize the relevant notions concerning directed hypergraphs (for more on directed hypergraphs, see \cite{MR1915421,MR1217096}). 
For us, a directed hypergraph $\mathcal{H}$ will be a pair $(V,\mathcal{E})$ where $V$ is a set of nodes and $\mathcal{E}$, the set of hyperedges, consists of pairs of multisets of nodes. An F-hyperedge is a hyperedge $e=(T(e),H(E))$ such that $|T(e)|=1$. An F-hypergraph is a directed hypergraph whose edges are F-hyperedges. 

\begin{definition} \label{hyper2} \normalfont
Given a tensor $A$ of order $m$ and rank $n$, we associate to it a weighted F-hypergraph 
\[
\mathcal{E}(A):=\{(i_1,\{i_2,i_3,\ldots,i_m\})\, : \, i_j \in [n] \, \forall j=1,2,\ldots,m\}
\]
whose vertices are $V(A):=[n]$, hyperarcs are $\mathcal{E}(A):=\{(i_1,\{i_2,i_3,\ldots,i_m\})\, : \, i_j \in [n] \, \forall j=1,2,\ldots,m\}$ and weights are given by $W(A)((i_1,\{i_2,\ldots,i_m\})=a_{i_1 i_2 \cdots i_m}$.
\end{definition}

Our study of H-eigenpairs motivates us to define the following type of cycle for an F-hypergraph whose hyperarc heads are of the same size.

\begin{definition} \label{HcycleDef} \normalfont
Let $\mathcal{H}=(V,\mathcal{E})$ be an F-hypergraph.
An \emph{H-cycle} in $\mathcal{H}$ is sequence 
\[
\{(v_{i_1^{(1)}},\{v_{i_2^{(1)}},\ldots,v_{i_m^{(1)}}\}),(v_{i_1^{(2)}},\{v_{i_2^{(2)}},\ldots,v_{i_m^{(2)}}\}),\ldots,(v_{i_1^{(r)}},\{v_{i_2^{(r)}},\ldots,v_{i_m^{(r)}}\})\}
\]
of hyperedges, such that the following combinatorial condition holds:
\[
\sum_{j=1}^r (m-1)e_{i_1^{(j)}}-e_{i_2^{(j)}}-\ldots-e_{i_m^{(j)}}=0.
\]
\end{definition}

We note that any tight cycle is also an H-cycle.

In this language, $\lambda(A)$ is equal to the minimal normalized weighted length of any H-cycle of $\mathcal{H}(A)$.

Finally, we introduce and study the H-cycle polytope, which is the feasible set of the linear program (\ref{LP}).

\begin{definition} \label{hcyclePolytope} \normalfont The \emph{H-cycle Polytope} $H_{n,m}$ is the polytope in $\R^{n^m}$ defined by the inequalities
\[
\left\lbrace
\begin{array}{llll}
\sum_{(i_1,i_2,\ldots,i_m) \in [n]^m} y_{i_1 i_2 \cdots i_m}((m-1)e_{i_1}-e_{i_2}-\ldots-e_{i_m})=0 \\
\sum_{(i_1,i_2,\ldots,i_m) \in [n]^m} y_{i_1 i_2 \cdots i_m}=1 \\
y_{i_1 i_2 \cdots i_m} \geq 0.
\end{array}
\right.
\]
The \emph{Fundamental H-cycles} are the vertices of $H_{n,m}$.
\end{definition}

The H-cycle polytope is a generalization of the normalized cycle polytope \cite{MR3201242} but turns out to be more complicated for $m>2$. For instance, the vertices are no longer normalized characteristic functions of subsets of the edges (see Example \ref{complicatedVertex}). Nevertheless, we show that it has the following nice property, which explains Theorem \ref{LPthm}.

\begin{theorem}\label{verticesOfHPolytope}
Let $y \in \R^{n^m}$ be a vertex of the H-cycle polytope $H_{n,m}$. Then $y$ has at most one nonzero entry of the form $y_{j i_2 \cdots i_m}$ for each $j=1,2,\ldots,n$. There exist vertices with $n$ nonzero entries. 
\end{theorem}

\section{H-eigenvalues of Tensors \label{eigenSection}}

 Our setup is that of Theorem \ref{uniquenessExistenceWithInfinity}:

\thmrestatement{\ref{uniquenessExistenceWithInfinity}}
Let $A \in \overline{\R}^{n^m}$ and assume that for each $i=1,2,\ldots,n$, the sets
\[
S_i:=\{\{i_2^{(i)},\ldots,i_m^{(i)}\} \, : \, a_{i i_2^{(i)} \cdots i_m^{(i)}} \neq \infty\}
\]
are nonempty and mutually equal. Then $A$ has a unique tropical H-eigenvalue $\lambda(A) \in \R$.

We begin by proving uniqueness of the H-eigenvalue. The main tool for the proof is to use Gordan's Theorem \cite[Theorem 2.2.1]{MR2184742} which states that for a matrix $M$
\[
\text{either} \quad
\exists x \in \mathbb{R}_+^n\setminus\{0\} \text{ such that } Mx = 0,
\quad\text{or}\quad
\exists y\in\mathbb{R}^n \text{ such that } M^t y > 0.
\]
The main idea in the proof is to perform Gaussian elimination using only additions of positive scalar multiples of rows.
\begin{proposition} \label{valueofEigenvalue}
Let $\lambda$ be a tropical H-eigenvalue of the tensor $A=(a_{i_1 \cdots i_m})$. For all $(i_1,i_2,\ldots,i_m) \in [n]^m$, define the column vector
\[
v_{i_1 \cdots i_m}:=-(m-1)e_{i_1}+e_{i_2}+\ldots+e_{i_m}
\] 
and let $M$ be the matrix whose columns are $v_{i_1 \cdots i_m}$. Then
\begin{equation} \label{LP1}
\lambda=\left[\begin{array}{ll@{}ll}
\text{minimize}  & \sum_{i_1 \cdots i_m}  a_{i_1 \cdots i_m} c_{i_1 \cdots i_m} &\\
\text{subject to}& Mc=0\\
& \mathbf{1}^t c=1 \\
& c \geq 0.
\end{array}\right]
\end{equation}
\end{proposition}

\proof
Let $x=(x_1,\ldots,x_n)^t$. If $\lambda$ is an H-eigenvalue, then we can write 
\[
\lambda=a_{i i_2^{(i)} \cdots i_m^{(i)}}+x^t v_{i i_2^{(i)} \cdots i_m^{(i)}}, \quad i=1,2,\ldots,n,
\]
for $a_{i i_2^{(i)} \cdots i_m^{(i)}} \in \R$.
We show that there exists $c \geq 0$, $c \neq 0$, such that
\begin{align} \label{positivec}
\sum_{i=1}^n c_{i} v_{i i_2^{(i)} \cdots i_m^{(i)}}=0.
\end{align}
Form the matrix $W$ whose columns are $v_{i i_2^{(i)} \cdots i_m^{(i)}}$. We would like to show that there exists $c \in \R^n_+ \setminus \{0\}$ such that $Wc=0$. Applying Gordan's Theorem, we will show that the alternative $\exists y\in\mathbb{R}^n \text{ such that } W^t y > 0$ leads to a contradiction. We do so as follows. If $W^t$ has a zero row, then we are done. So assume otherwise. We show that we can bring the matrix $W^t$ to row echelon form using the operation of adding to a row a positive scalar multiple of another row. This preserves the positivity condition. Since $W^t$ is not of full rank (the sum of columns adds to zero), this will prove the claim.

Observe that an entry of $W^t$ is negative if and only if it lies along the diagonal. Let $r_i=(r_{i1},r_{i2},\ldots,r_{in}), i=1,2,\ldots,n$ denote the rows of $V^t$. We may add positive scalar multiples of $r_1$ to all the remaining rows so as to zero out $r_{i1}$ for each $i \geq 2$. In doing so, either: the $(2,2)$ entry of $W^t$ remains negative, or the second row of $W^t$ is now zero. To see why this is so, note that the sum of entries of $r_2$ must be zero, while the entries $(2,3),\ldots,(2,n)$ are nonnegative. We proceed in this manner, adding a positive multiple of $r_i$ to $r_{i+1}$ until $W^t$ has reached row echelon form. The existence of $c \geq 0, c \neq 0$ satisfying  (\ref{positivec}) now follows.

 We then have
\[
\sum_{i=1}^n c_i \lambda= \sum_{i=1}^n c_i a_{i i_2^{(i)} \cdots i_m^{(i)}} +\sum_{i=1}^n c_i x^t v_{i i_2^{(i)} \cdots i_m^{(i)}}=\sum_{i=1}^n c_i a_{i i_2^{(i)} \cdots i_m^{(i)}},
\]
so that 
\[
\lambda \geq \left[\begin{array}{ll@{}ll}
\text{minimize}  & \sum_{i_1 \cdots i_m}  a_{i_1 \cdots i_m} c_{i_1 \cdots i_m} &\\
\text{subject to}& Mc=0\\
& \mathbf{1}^t c=1 \\
& c \geq 0
\end{array}\right].
\]

Conversely, from (\ref{eigenProblemH}), we have
\[
\lambda \leq a_{i_1 \cdots i_m}+x^tv_{i_1 \cdots i_m} 
\]
for every choice of indices. Suppose that $c \in \ker M, c \geq 0$ and $\mathbf{1}^t c =1$. Then, in particular,
\[
\sum_{i_1 \cdots i_m} c_{i_1 \cdots i_m} v_{i_1 \cdots i_m}=0.
\]
Since $c \geq 0$,
\[
(\sum_{i_1 \cdots i_m} c_{i_1 \cdots i_m}) \lambda(A) \leq \sum_{i_1 \cdots i_m} c_{i_1 \cdots i_m} a_{i_1 \cdots i_m}+\sum_{i_1 \cdots i_m} c_{i_1 \cdots i_m} x^t v_{i_1 \cdots i_m}=\sum_{i_1 \cdots i_m} c_{i_1 \cdots i_m} a_{i_1 \cdots i_m}.
\]
It follows that
\[
\lambda \leq \min_{\substack{c \in \ker M \\ c \geq 0, \mathbf{1}^t c =1}} \sum_{i_1 \cdots i_m} c_{i_1 \cdots i_m} a_{i_1 \cdots i_m}.
\]

\proofend

Next we show the existence of an H-eigenpair. The idea of the proof is to consider a tropical analogue of the proof of the Perron-Frobenius Theorem. For a reference on the Perron-Frobenius Theory of tensors, see \cite{MR2435198}.

For the proof of existence, we will be needing the following lemma, which implies that in tropical arithmetic, the ratio of two polynomials with the same support is bounded.

\begin{lemma} \label{polyRatioBounded}
Let $c_i,d_i,A_i \in \R$ for $1 \leq i \leq k$. Then 
\[
\min_i c_i-d_i \leq \min_i \{c_i+A_i\} - 
\min_i \{d_i+A_i\} \leq \max_i c_i-d_i .
\]
\end{lemma}

\proof
Suppose that 
\begin{align} \label{minEq}
\min_i \{c_i+A_i\}=c_i+p_i(x), \quad \min_i \{d_i+A_i\}=d_j+p_j(x).
\end{align}
By (\ref{minEq}), $c_i+A_i \leq c_j+A_j$ and $d_j+A_j \leq d_i+A_i$. Therefore
\[
\min_i \{c_i+A_i\}-\min_i \{d_i+A_i\}=c_i-d_j+A_i-A_j \leq    c_j-d_j,
\]
and similarly $c_i-d_j+A_i-A_j \geq c_i-d_i$.
\proofend

\begin{proposition}\label{existence}
For any tensor $A$ satisfying the hypotheses of Theorem \ref{uniquenessExistenceWithInfinity}, there exists a tropical H-eigenpair.
\end{proposition}

\proof 
Let $S:=S_i$. Define $F: \mathbb{T} \mathbb{P}^{n-1} \rightarrow \mathbb{T} \mathbb{P}^{n-1}$ by
\[
F(x)_i=\frac{\bigoplus_{\{i_2,\ldots,i_m\} \in S} a_{i i_2 \cdots i_m} \odot x_{i_2} \odot \cdots \odot x_{i_m}}{m-1}.
\]
Using the equivalence relation on $\mathbb{T} \mathbb{P}^{n-1}$, we can view this map as a map $F:\R^{n-1} \rightarrow \R^{n-1}$ by 
\[
F(x_1,x_2,\ldots,x_{n-1})_i=\frac{\bigoplus_{\{i_2,\ldots,i_m\} \in S} a_{i i_2 \cdots i_m} \odot x_{i_2} \odot \cdots \odot x_{i_m}-\bigoplus_{\{i_2,\ldots,i_m\} \in S} a_{n i_2 \cdots i_m} \odot x_{i_2} \odot \cdots \odot x_{i_m}}{m-1}
\]
for $1 \leq i \leq n-1$. By Lemma \ref{polyRatioBounded}, each coordinate of $F$ is bounded. Thus $F$ is a continuous mapping of a convex set of $\R^{n-1}$ into a bounded closed subset of $\R^{n-1}$, and consequently has a fixed point. This condition translates to having, for each $1 \leq i \leq n-1$,
\[
\bigoplus_{\{i_2,\ldots,i_m\} \in S} a_{i i_2 \cdots i_m} \odot x_{i_2} \odot \cdots \odot x_{i_m}=\bigoplus_{\{i_2,\ldots,i_m\} \in S} a_{n i_2 \cdots i_m} \odot x_{i_2} \odot \cdots \odot x_{i_m}+(m-1)x_i,
\]
i.e., the existence of an H-eigenpair (here $x_n$ is normalized to zero).
\proofend

\begin{remark}
The proof of existence can be easily adopted to tropical E-eigenpairs; one simply replaces the function 
\[
F(x)_i=\frac{\bigoplus_{\{i_2,\ldots,i_m\} \in S} a_{i i_2 \cdots i_m} \odot x_{i_2} \odot \cdots \odot x_{i_m}}{m-1}
\]
with
\[
F(x)_i=\bigoplus_{\{i_2,\ldots,i_m\} \in S} a_{i i_2 \cdots i_m} \odot x_{i_2} \odot \cdots \odot x_{i_m}.
\]
\end{remark}
Next we show that the H-eigenvalue of a tensor can be computed via a linear program.

\thmrestatement{\ref{LPthm}}
Assume the hypotheses of Theorem \ref{uniquenessExistenceWithInfinity}. The solution to the linear program
\begin{equation} \label{LP2}
\begin{array}{ll@{}ll}
\text{maximize}  & \lambda &\\
\text{subject to}& a_{i_1 i_2 \cdots i_m}+x_{i_2}+\ldots+x_{i_m} \geq \lambda + (m-1)x_{i_1}, &   & \forall (i_1,i_2,\ldots,i_m) \in [n]^m\\
\end{array}
\end{equation}
  is equal to the H-eigenvalue of $A$. Dually, the H-eigenvalue of $A$ is given by
  \begin{equation} \label{LPdual2}
\begin{array}{ll@{}ll}
\text{minimize}  & \sum_{(i_1,i_2,\ldots,i_m) \in [n]^m} a_{i_1 i_2 \cdots i_m} y_{i_1 i_2 \cdots i_m} &\\
\text{subject to}& \sum_{(i_1,i_2,\ldots,i_m) \in [n]^m} y_{i_1 i_2 \cdots i_m}((m-1)e_{i_1}-e_{i_2}-\ldots-e_{i_m})=0 \\
& \sum_{(i_1,i_2,\ldots,i_m) \in [n]^m} y_{i_1 i_2 \cdots i_m}=1 \\
& y_{i_1 i_2 \cdots i_m} \geq 0.
\end{array}
\end{equation}

\proof
The dual problem to (\ref{LP2}) is
\begin{equation} 
\begin{array}{ll@{}ll}
\text{minimize}  & \sum_{(i_1,i_2,\ldots,i_m) \in [n]^m} a_{i_1 i_2 \cdots i_m} y_{i_1 i_2 \cdots i_m} &\\
\text{subject to}& \sum_{(i_1,i_2,\ldots,i_m) \in [n]^m} y_{i_1 i_2 \cdots i_m}((m-1)e_{i_1}-e_{i_2}-\ldots-e_{i_m})=0 \\
& \sum_{(i_1,i_2,\ldots,i_m) \in [n]^m} y_{i_1 i_2 \cdots i_m}=1 \\
& y_{i_1 i_2 \cdots i_m} \geq 0.
\end{array}
\end{equation}
The primal problem is feasible, since we can take $x_i=0 \, \forall i \in [n]$ and $\lambda = \min_{(i_1,i_2,\ldots,i_m) \in [n]^m} a_{i_1 i_2 \cdots i_m}$. The problems (\ref{LP1}) and (\ref{LPdual2}) are identical, hence feasible, and so  the result follows by Proposition \ref{valueofEigenvalue}.
\proofend

For matrices ($m=2$), the dual problem has the interpretation of giving the minimal normalized length of any directed cycle in the weighted directed graph associated to the matrix $A$. Motivated by this, we associate to a tensor $A$ a weighted directed  hypergraph. Now follows a brief summary of the relevant notions concerning directed hypergraphs (for more on directed hypergraphs, see \cite{MR1915421,MR1217096}).

For us, a directed hypergraph $\mathcal{H}$ will be a pair $(V,\mathcal{E})$ where $V$ is a set of nodes and $\mathcal{E}$, the set of hyperedges, consists of pairs of multisets of nodes. An F-hyperedge is a hyperedge $e=(T(e),H(E))$ such that $|T(e)|=1$. An F-hypergraph is a directed hypergraph whose edges are F-hyperedges. 

\defrestatement{\ref{hyper2}} \normalfont
Given a tensor $A$ of order $m$ and rank $n$, we associate to it a weighted F-hypergraph $\mathcal{H}(A)=(V(A),\mathcal{E}(A),W(A))$ whose vertices are $V(A):=[n]$, hyperarcs are 
\[
\mathcal{E}(A):=\{(i_1,\{i_2,i_3,\ldots,i_m\})\, : \, i_j \in [n] \, \forall j=1,2,\ldots,m\}
\]
 and weights are given by $W(A)((i_1,\{i_2,\ldots,i_m\})=a_{i_1 i_2 \cdots i_m}$.

The linear program (\ref{LPdual}) motivates us to define the following type of cycle.
 
\defrestatement{\ref{HcycleDef}}
\normalfont
Let $\mathcal{H}=(V,\mathcal{E})$ be an F-hypergraph.
An \emph{H-cycle} in $\mathcal{H}$ is sequence 
\[
\{(v_{i_1^{(1)}},v_{i_2^{(1)}},\ldots,v_{i_m^{(1)}}),(v_{i_1^{(2)}},v_{i_2^{(2)}},\ldots,v_{i_m^{(2)}}),\ldots,(v_{i_1^{(r)}},v_{i_2^{(r)}},\ldots,v_{i_m^{(r)}})\}
\]
of hyperedges, such that the following combinatorial condition holds:
\[
\sum_{j=1}^r (m-1)e_{i_1^{(j)}}-e_{i_2^{(j)}}-\ldots-e_{i_m^{(j)}}=0.
\]

Recall that a tight cycle in a hypergraph is a sequence of edges of the form
\[
\{(v_1,v_2,\ldots,v_k),(v_2,v_3,\ldots,v_{k+1}),\ldots,(v_{r-k+1},\ldots,v_{r-1},v_r),(v_{r-k+2},\ldots,v_r,v_1),\ldots,(v_{r},v_1,\ldots,v_{k-1})\},
\]
with the vertices not necessarily distinct. Any tight cycle is also an H-cycle. Indeed, considering indices modulo $r$,
\[
\sum_{j=1}^r (m-1)e_{v_j}-e_{v_{j+1}}-\ldots-e_{v_{j+m-1}}=(m-1)(e_{v_1}+\ldots+e_{v_r})
 \underbrace{-(e_{v_1}+\ldots+e_{v_r})-\ldots-(e_{v_1}+\ldots+e_{v_r})}_\text{m-1 times} =0.
\]

In this language, (\ref{LPdual}) shows that $\lambda(A)$ is equal to the minimal normalized weighted length of any H-cycle of $\mathcal{H}(A)$.

\section{H-eigenvalues of Symmetric Tensors}

In this section, we specialize the results of section \ref{eigenSection} to symmetric tensors.

\begin{proposition}\label{symmetricEVal}
For a symmetric tensor $A \in \R^{n^m}$,
\[
\lambda(A)=\min_{i_1 i_2 \cdots i_m} a_{i_1 i_2 \cdots i_m}.
\]
\end{proposition}

\proof

Consider the system (\ref{LP}) for an H-eigenpair. Set $\Delta_{i,j}=x_i-x_j$. In this notation, we have
\[
a_{i_1 i_2 \cdots i_m}+\Delta_{i_2,i_1}+\Delta_{i_3,i_1}+\ldots+\Delta_{i_m,i_1} \geq \lambda(A)  \quad \forall i_1,\ldots,i_m \in [n].
\]
In particular, for every $\sigma \in S_m$, we have
\[
\lambda(A) \leq a_{i_{\sigma(1)} i_{\sigma(2)} \cdots i_{\sigma(m)}}+\Delta_{i_{\sigma(2)},i_{\sigma(1)}}+\Delta_{i_{\sigma(3)},i_{\sigma(1)}}+\ldots+\Delta_{i_{\sigma(m)},i_{\sigma(1)}}.
\]
Summing such inequalities over $S_m$ and using the symmetry of $A$ shows that $\lambda(A) \leq \min_{i_1 i_2 \cdots i_m} a_{i_1 i_2 \cdots i_m}$. Taking $x_i=0 \, \forall i=1,2,\ldots,n$ in (\ref{LP}) and $\lambda=\min_{i_1 i_2 \cdots i_m} a_{i_1 i_2 \cdots i_m}$ yields the result.
\proofend

\begin{proposition}
Suppose that $A$ is a symmetric tensor with  minimum entries $a_{i_1 i_2 \cdots i_m}$ having the same set of indices $I=\{i_1,i_2,\ldots,i_m\}$.  If $x$ is an eigenvector of $A$, then
\begin{enumerate}
\item \label{equalI} $x_i=x_j$ for each $i,j \in I$.
\item \label{lessI} $x_i \leq x_j$ for each $i \in I$ and $j \not \in I$.
\end{enumerate}
\end{proposition}

\proof
Let $x_k= \min_l x_l$. For some choice of indices $k_2,\ldots,k_m$,
\[
a_{k k_2 \cdots k_m}+\Delta_{k_2,k}+\Delta_{k_3,k}+\ldots+\Delta_{k_m,k} = \lambda.
\]
 By choice of $k$, $\Delta_{k,k_l} \geq 0$ for $l=2,3,\ldots,m$. By Proposition \ref{symmetricEVal}, we must have $\Delta_{k,k_l}=0$ for each $l$. Since $\lambda(A)$ is equal to the minimum entry of $A$, $a_{k k_2 \cdots k_m}$ is a minimum element. Since a minimum element $a_{j_1 \cdots j_m}$ has indices $I$, $\{k,k_2,\ldots,k_m\}=I$. This shows (\ref{equalI}). Part (\ref{lessI}) follows from $k \in I$, $x_k$ being minimal and part (\ref{equalI}).
\proofend

\section{H-cycle Polytope}

\defrestatement{\ref{hcyclePolytope}} \normalfont The \emph{H-cycle Polytope} $H_{n,m}$ is the polytope in $\R^{n^m}$ defined by the inequalities
\[
\left\lbrace
\begin{array}{llll}
\sum_{(i_1,i_2,\ldots,i_m) \in [n]^m} y_{i_1 i_2 \cdots i_m}((m-1)e_{i_1}-e_{i_2}-\ldots-e_{i_m})=0 \\
\sum_{(i_1,i_2,\ldots,i_m) \in [n]^m} y_{i_1 i_2 \cdots i_m}=1 \\
y_{i_1 i_2 \cdots i_m} \geq 0.
\end{array}
\right.
\]

In view of (\ref{LPdual}), we can interpret tensors as linear functionals on the H-cycle polytope. 

The H-cycle polytope is a generalization of the normalized cycle polytope \cite{MR3201242} but turns out to be more complicated for $m>2$. For instance, the vertices are no longer normalized characteristic functions of subsets of the edges, as the following example shows. Nevertheless, we will show that it has some nice properties.

\begin{example} \label{complicatedVertex}
Consider the tensor whose entries are all zero except for $a_{132}=a_{213}=a_{322}=-1$. Taking
\begin{align}\label{vertex1} 
y_{132}=\frac{2}{9}, \quad y_{213}=\frac{4}{9}, \quad y_{322}=\frac{3}{9}, \quad y_{ijk}=0 \text{ for all remaining indices}
\end{align}
we see that the dual problem (\ref{LPdual}) attains a value of $-1$. The primal problem (\ref{LP}) also has a feasible point attaining the value $-1$: we take $\lambda=-1$ and $x_i=0$ for each $i$. It is easy to check that (\ref{vertex1}) is a vertex (not only a face) of the H-cycle polytope.
\end{example}

The following result shows that the coordinate vectors of the vertices of the H-cycle polytope have some nice structure.

\thmrestatement{\ref{verticesOfHPolytope}} Let $y \in \R^{n^m}$ be a vertex of the H-cycle polytope $H_{n,m}$. Then $y$ has at most one nonzero entry of the form $y_{j i_2 \cdots i_m}$ for each $j=1,2,\ldots,n$. There exist vertices with $n$ nonzero entries.

\proof
Let $y$ be a vertex of $H_{n,m}$. Take a linear functional $A=(a_{i_1 i_2 \cdots a_m})$ whose restriction to $H_{n,m}$ is minimized on, and only on, $y$. This minimum value is equal to $\lambda^*$, the solution to the linear program (\ref{LPdual}). By Theorem \ref{LPthm}, the value $\lambda(A)$ is attained by a tropical H-eigenvector of the tensor $A$. From the proof of Proposition \ref{valueofEigenvalue}, we see that the tropical H-eigenvector $\lambda$ can be written as 
\[
\lambda= \sum_{i=1}^n c_i a_{i i_2^{(i)} \cdots i_m^{(i)}}, \quad c \in H_{n,m}.
\]
It follows that $y_{i i_2^{(i)} \cdots i_m^{(i)}}=c_i$ and the remaining entries of $y$ are zero.

Next we show that a vertex with $n$ nonzero entries exists. Consider the (one-sided) infinite sequence
\[
s=\{\overline{1,2,\ldots,n}\}.
\]
We take $y$ to have an entry of $\frac{1}{n}$ at position $(s_j, s_{j+1},\ldots, s_{j+m})$ for $j=1,2,\ldots,n$, and zero otherwise.
We also take the linear functional $A$ to have entries of $-1$ at these positions and zeros otherwise. To see that $y \in H_{n,m}$, write $s_j=j-n \lfloor \frac{j-1}{n} \rfloor$.
Then 
\[
\sum_{(i_1,i_2,\ldots,i_m) \in [n]^m} y_{i_1 i_2 \cdots i_m}((m-1)e_{i_1}-e_{i_2}-\ldots-e_{i_m})
\]
\[
=\sum_{j=1}^n \frac{1}{n} ((m-1) e_{j-n \lfloor \frac{j-1}{n} \rfloor }-e_{j+1-n \lfloor \frac{j}{n}\rfloor}-\ldots-e_{j+m-n \lfloor \frac{j+m-1}{n} \rfloor})
\]
\[
= \frac{1}{n} ((m-1)(e_1+\ldots+e_n)-(e_1+\ldots+e_n)-\ldots-(e_1+\ldots+e_n))=0.
\]
The dual problem (\ref{LPdual}) takes on the value $-1$. The primal problem (\ref{LP}) has a solution $\lambda=-1$, $x_i=0 \, \forall i$. This shows that $y$ lies on a face of $H_{n,m}$. To see that it is a vertex, note that by the rearrangement inequality, any other minimizer must have support contained in the support of $y$. We show that the matrix $W$ whose columns are $(m-1) e_{j-n \lfloor \frac{j-1}{n} \rfloor }-e_{j+1-n \lfloor \frac{j}{n}\rfloor}-\ldots-e_{j+m-n \lfloor \frac{j+m-1}{n} \rfloor}, j=1,2,\ldots,n$, has rank $n-1$. This implies that the kernel is $1$-dimensional, so that $y$ is indeed a vertex. Write $m=an+r, 1 \leq r \leq n$. The matrix $W$ is the circulant matrix whose first row is $(a-(m-1))e_1^t+(a+1)e_2^ t+\ldots+(a+1) e_r^ t+a e_{r+1}^t+\ldots+a e_n^t)$. The associated polynomial is
\[
f(x)=a-(m-1)+(a+1)x+\ldots+(a+1)x^{r-1}+a x^{r}+\ldots+a x^n
\]
\[
=a-(m-1)+x+\ldots+x^{r-1}+a(x+\ldots+x^{n-1}).
\]
It is known that the rank of a circulant matrix is equal to $n-d$, where $d$ is the degree of $\gcd(f(x),x^n-1)$ \cite{MR0080623}. We show that the only common root to both polynomials is $x=1$. We have
\[
f(1)=a-(m-1)+r-1+a(n-1)=0.
\]
 Let $\xi \neq 1$ be an $n$th root of unity. Then
 \[
 f(\xi)=-m+1+\xi+\ldots+\xi^{r-1}=\frac{1-\xi^r}{1-\xi}-m.
 \]
We see that $f(\xi) \neq 0$ for $m \geq 2$.
\proofend

\begin{corollary}
Let $s_j=j-n \lfloor \frac{j-1}{n} \rfloor$. The cycle whose edges are
\[
\{(s_1, s_{2},\ldots, s_{1+m}),(s_2, s_{3},\ldots, s_{2+m}),\ldots,(s_n, s_{n+1},\ldots, s_{n+m})\}
\]
is a vertex of the H-cycle polytope $H_{n,m}$.
\end{corollary}

\proof
This is a consequence of the proof of Theorem \ref{verticesOfHPolytope}.
\proofend

\section{Future Directions}

One direction deserving of study is the theory of the H-eigenvectors of a tensor. For non-generic tensors, there may be infinitely many nonequivalent H-eigenvectors. However, for generic tensors, one expects a unique H-eigenpair. 

Another direction of study is properties of the H-cycle polytope. For instance, as a start it would be of interest to understand the set of vertices better: its cardinality and characteristics. Next one may study the facial structure more generally.

It is also of interest to understand the relationship between the classical H-eigenpairs of a tensor and the tropical ones. For instance, to which classical H-eigenpairs do tropical ones lift.

\bigskip
{\bf Acknowledgments}.  
The author is grateful to Bernd Sturmfels for suggesting this project and providing valuable advice throughout. The author is also grateful to Jean-Christophe Yoccoz for helpful suggestions.

This material is based upon work supported by the National Science Foundation Graduate Research Fellowship under Grant No. DGE 1106400. Any opinion, findings, and conclusions or recommendations expressed in this material are those of the authors(s) and do not necessarily reflect the views of the National Science Foundation.

\bibliographystyle{alpha}
\bibliography{bibliography}

\end{document}